\documentclass[draft]{amsart}

\usepackage{amssymb}
\usepackage{url}
\usepackage{amscd}
\usepackage[all]{xy}

\theoremstyle{plain}
\newtheorem{theorem}{Theorem}

\newtheorem{corollary}[theorem]{Corollary}
\newtheorem{lemma}[theorem]{Lemma}

\newtheorem{definition}[theorem]{Definition}
\theoremstyle{remark}
\newtheorem{remark}{Remark}
\newtheorem{example}{Example}

\begin{document}

\date{}

\title{$*-$open sets and $*-$ continuity in topological spaces}

\author{Aliakbar Alijani}

\address{Mollasadra Technical and Vocational College\\
Technical and Vocational University\\
Ramsar, Iran}

\email{alijanialiakbar@gmail.com}

\thanks{}


\begin{abstract}
In this paper, we study some properties of $*-$open and $*-$closed subsets of a space. The collection of all $*-$open subsets of a space $X$ form a topology on $X$ which is denoted by $^{*}O(X)$. We investigate the relations between topological properties of $X$ with the topology $^{*}O(X)$ and $X$. Also, we introduce the concept of a $*-$continuous map.
\end{abstract}

\maketitle


\newcommand\sfrac[2]{{#1/#2}}

\newcommand\cont{\operatorname{cont}}
\newcommand\diff{\operatorname{diff}}


\section{Introduction}
\newcommand{\stk}{\stackrel}
By a space, we mean a topological space. In a space $X$, the closure and the interior of a set $A$ will be denoted by $cl(A)$ and $int(A)$, respectively. If $A\subseteq Y\subseteq X$, the closure and the interior of $A$ as a subspace of $Y$ will be denoted by $cl_{Y}(A)$ and $int_{Y}(A)$, respectively. A set $A$ in $X$ is called regular open if $int(cl(A))=A$. The collection of all regular open subsets of a space $X$ denote by $R.O(X)$. A set $A$ in space $X$ is called $*-$open if for every $x\in A$, there exists $ U\in R.O(X)$ containing $x$ such that $U\subseteq A$. We denote by $^{*}O(X)$, the collection of all $*-$open subsets of a space $X$. Then, $^{*}O(X)$ is a topology on $X$ (Theorem \ref{32}).  Recall that $R.O(X)$, form a base for a topology on $X$ \cite[p. 138]{B}. In fact, $^{*}O(X)$, is the topology generated by $R.O(X)$. In Section \ref{sec2}, We investigate the relations between some topological properties of $(X,^{*}O(X))$ and $X$ as follows:
\begin{itemize}
\item If $X$ is a compact space, then $(X,^{*}O(X))$ is a compact space (Lemma \ref{19}). Conversely is not true (Example \ref{20}).
\item If $X$ is a Hausdorrf, locally compact space, then $(X,^{*}O(X))$ is a Hausdorrf, locally compact space (Lemma \ref{21}). Conversely is not true (Example \ref{22}).
\item $X$ is a connected space if and only if $(X,^{*}O(X))$ is a connected space (Lemma \ref{18}).
\item If $X$ is a path connected space, then $(X,^{*}O(X))$ is a path connected space (Lemma \ref{23}). Conversely is not true (Example \ref{24}).
\item If $X$ is a locally connected space, then $(X,^{*}O(X))$ is a locally connected space (Lemma \ref{25}). Conversely is not true (Example \ref{27}).
\end{itemize}

In Section \ref{sec3}, we introduce the concept of $*-$continuous maps. A map $f:X\to Y$ is called $*-$continuous if the inverse image of every $*-$open set is $*-$open. $*-$continuous neither implies continuous nor is implied by continuous (Example \ref{28}). $*-$continuous image of connected sets are connected (Lemma \ref{29}).  compactness, path connectedness and local connectedness need not to be preserved by $*-$continuous maps (Example \ref{30}).

For two spaces $X$ and $Y$, We denote by $X\times Y$, the cartesian product of $X\times Y$ with the product topology. For more information on topological spaces, see \cite{B}.


\section{$*-$open and $*-$closed sets}\label{sec2}
In this section, we introduce the concept and study some properties of $*-$open and $*-$closed sets.

\begin{definition}
Let $X$ be a space. For a set $A$ in $X$, define $^{*}int(A)$ and $^{*}cl{A}$ as follows:

 $$^{*}int(A)=\{x\in A| \exists W\in R.O(X)  ; x\in W\subseteq A\}$$
  $$^{*}cl(A)=\{x\in X| \forall x\in W\in R.O(X); W\bigcap A\neq\emptyset\}$$
\end{definition}

\begin{definition}
A set $A$ in a space $X$ is called $*-open$ if $^{*}int(A)=A$ and $*-$closed if $^{*}cl{A}=A$.
\end{definition}

In Lemma \ref{1}, we present some properties of $^{*}int(A)$ and $^{*}cl(A)$.

\begin{lemma}\label{1}
Let $X$ be a  space and $A$ an arbitrary set in $X$. Then:
\begin{enumerate}
\item $^{*}int(A)\subseteq int(A)\subseteq A\subseteq cl(A)\subseteq ^{*}cl{A}$
\item $^{*}int(A)=\bigcup\{ W\in R.O(X); W\subseteq A\}.$
\item $^{*}cl(A)=\bigcap\{W\in R.C(X); A\subseteq W\}.$
\item $^{*}int(A)$ is open and $^{*}cl(A)$ closed.
\item $^{*}cl(X-A)=X-^{*}int(A)$ and $X-^{*}cl(A)=^{*}int(X-A)$.\\
\end{enumerate}
\end{lemma}

\proof
It is clear.
\endproof

\begin{remark}
It follows from (5) of Lemma \ref{1} that a set $A$ is $*-$open if and only if $X-A$ be  $*-$closed.
\end{remark}

\begin{remark}
 Clearly, $*-$open sets are open. But, conversely need not be true. See Example \ref{2}.
\end{remark}

\begin{example}\label{2}
 Let $\mathbb{R}$ be the set of reals and, $\tau$ be the topology on $\mathbb{R}$ generated by the union of $\tau_{1}$, the usual topology on $\mathbb{R}$, and $\tau_{2}$, the topology of countable complements on $\mathbb{R}$ \cite[Example 63]{SS}. Then, $\mathbb{R-Q}$ is open in $(\mathbb{R},\tau)$. But, $^{*}int(\mathbb{R-Q})=\emptyset$.\\
\end{example}

\begin{definition}
A space is called semi-regular if it has a basis consisting of regular open sets \cite{S}.
\end{definition}

\begin{lemma}\label{3}
Let $X$ be a semi-regular space. Then, a set $A$ in $X$ is $*-$open if and only if it is open.
\end{lemma}

\proof
Let $A$ be an open subset of $X$ and $x\in A$. Since $X$ is semi-regular, there exists a set $W\in R.O(X)$ such that $x\in W\subseteq A$. Hence, $x\in ^{*}Int(A)$ and $A$ is $*-$open.
\endproof

\begin{remark}
By (2) of Lemma \ref{1}, regular open sets are $*-$open. But, conversely need not be true. See Example \ref{4}.
\end{remark}

\begin{example}\label{4}
Since $(\mathbb{R},\tau_{1})$ is a semi-regular space, $(0,1)\bigcup(1,2)$ is $*-$open, but it is not regular open.
\end{example}

\begin{definition}
$^{*}O(X)$ will denote the class of all $*-$open sets in $X$.
\end{definition}

\begin{lemma}\label{5}
Let $(X,\tau)$ be a space. Then, $R.O(X)\subseteq ^{*}O(X)\subseteq \tau$.
\end{lemma}

\begin{lemma}\label{6}
Let $(X,\tau)$ be a semi-regular space. Then, $\tau=^{*}O(X)$.
\end{lemma}

\proof
It is clear by Lemma \ref{3}.
\endproof

\begin{remark}
In general, for a space $(X,\tau)$, $\tau\neq ^{*}O(X)$. It is sufficient to consider $(\mathbb{R},\tau)$ as in Example \ref{2}. Then, $^{*}O(X)=\tau_{1}\neq \tau$.
\end{remark}

\begin{lemma}\label{7}
Let $U$ be an open subset of $X$. Then, $Int(cl(U))$ is a regular open set.
\end{lemma}

\proof
It is clear.
\endproof

\begin{lemma}\label{8}
Let $(X,\tau)$ be a space. Then, $R.O(X,^{*}O(X))\subseteq R.O(X,\tau)$.
\end{lemma}

\proof
Let $W$ be a regular open set in $(X,^{*}O(X))$. Then,$ Int_{^{*}O(X)}(cl_{^{*}O(X)}(W))=W$. Let $x\in Int(cl(W))$. Then, $x\in U$ for some $U\in \tau$ and $U\subseteq cl(W)$. So, $x\in Int(cl(U))\subseteq cl(W)$. Since $^{*}O(X)\subseteq \tau$, $cl(W)\subseteq cl_{^{*}O(X)}(W)$. By Lemma \ref{7}, $Int(cl(U))\in ^{*}O(X)$. Hance, $x\in Int_{^{*}O(X)}(cl_{^{*}O(X)}(W))$. So, $x\in W$.
\endproof

\begin{lemma}
Let $(X,\tau)$ be a space. Then, $^{*}O(X,^{*}O(X,\tau))\subseteq^{*}O(X,\tau)$.
\end{lemma}

\proof
It is clear by Lemma \ref{8}.
\endproof

\begin{lemma}\label{9}
Let $Y$ be an open subset of $X$ and $A\subseteq Y$. Then we have:
\begin{enumerate}
\item $W\in R.O(Y)$ if and only if $W=N\bigcap Y$ where $N\in R.O(X)$.
\item $^{*}int(A)=^{*}int_{Y}(A)\bigcap ^{*}Int(Y)$
\item $^{*}cl_{Y}(A)=^{*}cl(A)\bigcap Y$
\end{enumerate}
\end{lemma}

\proof
(1) First suppose that $W\in R.O(Y)$. So, $int_{Y}(cl_{Y}(W))=W$. We have $$W=int_{Y}(cl_{Y}(W))=int(cl(U)\bigcap Y)=int(cl(U)\bigcap Y$$ Put $N=int(cl(U)$. Then, $N\in R.O(X)$ and the proof is complete. Conversely, let $W=N\bigcap Y$ where $N\in R.O(X)$. Then $$int_{Y}(cl_{Y}(W))=int(cl(N\bigcap Y)\bigcap Y)\subseteq int(cl(N))\bigcap Y=N\bigcap Y=W$$.
So, $W\in R.O(Y)$.
\endproof

\begin{corollary}\label{10}
Let $A\subseteq Y\subseteq X$ where $X$ is a space and $Y$ an open subspace. If $A\in ^{*}O(X)$, then $A\in ^{*}O(Y)$.
\end{corollary}

\proof
First recall that $^{*}Int(A)\subseteq ^{*}Int(Y)$. Since $A\in ^{*}O(X)$, $^{*}Int(A)=A$. Hence,$^{*}Int_{Y}(A)\subseteq A\subseteq ^{*}Int(Y)$. By (2) of Lemma \ref{9},$A\in ^{*}O(Y)$.
\endproof

\begin{remark}
The converse of Corollary \ref{10} need not be true. See Example \ref{11}.
\end{remark}

\begin{example}\label{11}
Let $(\mathbb{R},\tau)$ be as Example \ref{2}. Then, $\mathbb{R-Q}\in ^{*}O(\mathbb{R-Q})$. But, $\mathbb{R-Q}\notin ^{*}O(\mathbb{R-Q})$.
\end{example}



\begin{lemma}\label{12}
Let $X$ and $Y$ be two spaces,$A\in ^{*}O(X)$ and $B\in ^{*}O(Y)$. Then, $A\times B\in ^{*}O(X\times Y)$.
\end{lemma}

\proof
Let $A\in ^{*}O(X)$, $B\in ^{*}O(Y)$ and $(x,y)\in A\times B$. Then, there exists $x\in U\in R.O(X)$ and $y\in W\in R.O(Y)$ such that $U\subseteq A$ and $V\subseteq B$. Since $Int(cl(U\times V))=Int(cl(U))\times Int(cl(V))$, so $U\times V \in R.O(X\times Y)$. The proof is complete.
\endproof

\begin{remark}
It is not necessary that if $A\in ^{*}O(X\times Y)$, then $A=A_{1}\times A_{2}$ where $A_{1}\in ^{*}O(X)$ and $A_{2}\in ^{*}O(Y)$. See Example \ref{13}.
\end{remark}

\begin{example}\label{13}
Let $\mathbb{R}$ be the reals with the usual topology. By Lemma \ref{6}, $A=((0,1)\times (0,1))\cup (\mathbb{R}\times \{0\})$ is $*-$open in $\mathbb{R}\times \mathbb{R}_d$ ($\mathbb{R}_{d}$ is the reals with the discrete topology). An easy calculation shows that $A\neq A_{1}\times A_{2}$ for every $A_{1}\in ^{*}O(\mathbb{R})$ and $A_{2}\in ^{*}O(\mathbb{R}_d)$.
\end{example}

\begin{lemma}\label{14}
Let $\{A_{i};i\in I\}$ be a collection of $*-$open sets in a space $X$. Then, $\bigcup_{i}A_{i}$ is $*-$open set. If $I$ is finite, then $\bigcap_{i}A_{i}$ is $*-$open.
\end{lemma}

\proof
Let $x\in \bigcup_{i}A_{i}$. Then, $x\in A_{i}$ for some $i\in I$. So, there exists $U\in R.O(X)$ containing $x$ such that $U\subseteq A_{i}$. It is clear that $ U\subseteq \bigcup_{i}A_{i}$. Now, suppose that $x\in \bigcap_{i}A_{i}$. Then, $x\in A_{i}$ for each $i\in I$. So, there exists $U_{i}\in R.O(X)$ such that $U_{i}\subseteq A_{i}$ for each $i\in I$. Since $I$ is finite, $\bigcap_{i}U_{i}\in R.O(X)$ and $\bigcap_{i}U_{i}\subseteq \bigcap_{i}A_{i}$.
\endproof

\begin{remark}
The arbitrary intersection of $*-$open sets need not be a $*-$open set. See Example \ref{15}.
\end{remark}

\begin{example}\label{15}
In $(\mathbb{R},\tau_{1})$, it is clear that $\bigcap_{n\in \mathbb{N}}(-\frac{1}{n},\frac{1}{n})=\{0\}$ which is not $*-$open.
\end{example}

\begin{remark}
The converse of Lemma \ref{14} need not be true. See Example \ref{16}.
\end{remark}

\begin{example}\label{16}
Let $\mathbb{R}$ be the reals with the usual topology. Then, $(\mathbb{R-Q})\bigcup \mathbb{Q}=\mathbb{R}$ and $(\mathbb{R-Q})\bigcap \mathbb{Q}=\emptyset$ which are $*-$open sets. But, $\mathbb{R-Q}$ and $\mathbb{Q}$ are not $*-$open sets.
\end{example}

\begin{theorem}\label{32}
$^{*}O(X)$ is a topology on $X$.
\end{theorem}

\proof
Clearly, $\{\emptyset,X\}\subseteq ^{*}O(X)$. By Lemma \ref{14}, $^{*}O(X)$ is closed under finite intersections and arbitrary unions.
\endproof




\begin{remark}
$*-$open sets do not preserve by open continuous maps. See Example \ref{17}.
\end{remark}

\begin{example}\label{17}
Let $X$ be an infinite set and $p\in X$. Define $\tau=\{\phi\}\bigcup \{U\subseteq X;p\in U\}$. Then, $\tau$ is a topology on $X$. Consider $f:X\to X$ as follows: $f(x)=p$ for all $x\in X$. Then, $f$ is an open continuous map. But, $f(X)=\{p\}\notin ^{*}O(X)$.
\end{example}





\begin{lemma}\label{18}
Let $(X,\tau)$ be a space. Then, $(X,\tau)$ is a connected space if and only if $(X,^{*}O(X))$ be a connected space.
\end{lemma}

\proof
Let $(X,\tau)$ be connected. By Lemma \ref{5}, it is clear that $(X,^{*}O(X))$ is connected. Conversely, let $(X,^{*}O(X))$ be connected and $A$, a clopen subset of $(X,\tau)$. Then, $\{A, X-A\}\subseteq ^{*}O(X)$. So, $A$ is a clopen set in $(X,^{*}O(X))$. Since $(X,^{*}O(X))$ is connected, $A=\emptyset$ or $A=X$ and proof is complete.
\endproof

\begin{lemma}\label{19}
Let $(X,\tau)$ be compact space. Then, $(X,^{*}O(X))$ is a compact space.
\end{lemma}

\proof
It is clear by Lemma \ref{5}.
\endproof

\begin{remark}
The converse of Lemma \ref{19} need not be true. See Example \ref{20}.
\end{remark}

\begin{example}\label{20}
Let $X$ be as Example \ref{17}. Then, $^{*}O(X)=\{\emptyset,X\}$. So, $(X,^{*}O(X))$ is compact. But, $(X,\tau)$ is not compact.
\end{example}

\begin{lemma}\label{21}
Let $(X,\tau)$ be a Hausdorff, locally compact space. Then, $(X,^{*}O(X))$ is locally compact.
\end{lemma}

\proof
Let $x\in X$. Then, $U\subseteq C$ for some $x\in U\in \tau$ and compact set $C$ in $(X,\tau)$. By Lemma \ref{7}, $Int(cl(U))\in ^{*}O(X)$ and $Int(cl(U))\subseteq C$. By Lemma \ref{5}, it is clear that $C$ is compact in $(X,^{*}O(X))$.
\endproof

\begin{remark}
The converse of Lemma \ref{21} need not be true. See Example \ref{22}.
\end{remark}

\begin{example}\label{22}
Let $(\mathbb{R},\tau)$ be the same space as in Example \ref{2}. Then, $^{*}O(\mathbb{R})=\tau_{1}$. So, $(\mathbb{R},^{*}O(\mathbb{R}))$ is locally compact. But, $(\mathbb{R},\tau)$ is not locally compact.
\end{example}

\begin{lemma}\label{23}
Let $(X,\tau)$ be a path connected space. Then, $(X,^{*}O(X))$ is a path connected space.
\end{lemma}

\proof
Let $(X,\tau)$ be a path connected space and $x_{0},x_{1}\in X$. Then, there exists a continuous map $f:I\to X$ such that $f(0)=x_{0}$ and $f(1)=x_{1}$. Since the identity map $1_{X}:(X,\tau)\to (X,^{*}O(X))$ is continuous, $1_{X}of:I \to (X,^{*}O(X))$ is continuous. So, $(X,^{*}O(X))$ is a path connected.
\endproof

\begin{remark}
The converse of Lemma \ref{23} need not be true. See Example \ref{24}.
\end{remark}

\begin{example}\label{24}
Let $(\mathbb{R},\tau)$ be the same space as in Example \ref{2}. Since $^{*}O(\mathbb{R})=\tau_{1}$, $(\mathbb{R},^{*}O(\mathbb{R}))$ is path connected. Now suppose that, $x_{0},x_{1}\in \mathbb{R}$ and $f:I \to (\mathbb{R},\tau)$ be a continuous map such that $f(0)=x_{0}$ and $f(1)=x_{1}$. Then, $f(I)$ is a compact, connected subset of $(\mathbb{R},\tau)$ which is a contradiction.
\end{example}

\begin{lemma}\label{25}
Let $(X,\tau)$ be a locally connected space. Then, $(X,^{*}O(X))$ is locally connected.
\end{lemma}

\proof
Let $(X,\tau)$ be a locally connected space and $x\in U\in ^{*}O(X)$. Then, there exist $V\in R.O(X)$ containing $x$ such that $V\subseteq U$. Since $(X,\tau)$ is locally connected, there exists $ W\in \tau$ containing $x$ such that $W$ is connected and $W\subseteq V$. Since $W\subseteq Int(cl(W))\subseteq cl(W)$, $ Int(cl(W))$ is connected in $(X,\tau)$. It is clear that $Int(cl(W))$ is connected in $(X,^{*}O(X))$ and proof is complete.
\endproof

\begin{lemma}\label{26}
Let $(\mathbb{R},\tau)$ be as in Example \ref{2}. A set $A$ is connected in $(\mathbb{R},\tau)$ if and only if $A$ is connected in $(\mathbb{R},^{*}O(\mathbb{R}))$
\end{lemma}

\proof
Let $A$ be a set in $(\mathbb{R},\tau)$, $A=U_{1}\bigcup U_{2}$ where $U_{1},U_{2}\in \tau$ and $U_{1}\bigcap U_{2}=\emptyset$. Then, $U_{1},U_{2}$ are closed in $(\mathbb{R},\tau)$. Hence, $U_{1},U_{2}\in (\mathbb{R},\tau_{1})$ and proof is complete.
\endproof

\begin{remark}
The converse of Lemma \ref{25} need not be true. See Example \ref{27}.
\end{remark}

\begin{example}\label{27}
Let $(\mathbb{R},\tau)$ be the same space as in Example \ref{2}. Consider the open set $\mathbb{R-Q}$ in $(\mathbb{R},\tau)$. If $U$ be an open connected set in $(\mathbb{R},\tau)$ such that $U\subseteq \mathbb{R-Q}$, then by Lemma \ref{26}, $U\in \tau_{1}$ which is a contradiction (since $U\bigcap\mathbb{Q}\neq\emptyset$). So, $(\mathbb{R},\tau)$ is not locally connected. Recall that $^{*}O((\mathbb{R},\tau))=\tau_{1}$ and $(\mathbb{R},\tau_{1})$ is locally connected.
\end{example}

\section{$*-$continuous maps}\label{sec3}
In this section, we introduce the concept and study some properties of $*-$continuous maps.

\begin{definition}
A function $f:X\to Y$ is said to be $*-$continuous if the inverse image of every $*-$open subset of $Y$ is $*-$open subset of $X$.
\end{definition}

\begin{remark}
Recall that $*-$continuous neither implies continuous nor is implied by continuous. See Example \ref{28}.
\end{remark}

\begin{example}\label{28}
Let $(\mathbb{R},\tau)$ be as in Example \ref{2}. The identity map from $(\mathbb{R},\tau_{1})$ to $(\mathbb{R},\tau)$ is $*-$continuous. But, it is not continuous. Also, consider $f:(\mathbb{R},\tau)\to(\mathbb{R},\tau_{1})$ as follows: $f(x)=x$ if $x\in \mathbb{R-Q}$ and $f(x)=0$ if $x\in \mathbb{Q}$. Then, $f$ is continuous since if $a<b$ and $0\in (a,b)$, then $f^{-1}((a,b))=(a,b)$ and otherwise $f^{-1}((a,b))=(a,b)-\mathbb{Q}$. $f$ is not $*-$continuous, since $(0,1)$ is $*-$open in $(\mathbb{R},\tau_{1})$ but $f^{-1}((0,1))=(0,1)-\mathbb{Q}$ is not $*-$open in $(\mathbb{R},\tau)$.
\end{example}

\begin{lemma}
The projection maps $\pi_{i}:\prod_{j}X_{j}\to X_{i}$ are $*-$continuous.
\end{lemma}

\proof
For $i\in I$, let $U_{i}\in ^{*}O(X_{i})$. Then, $\pi_{i}^{-1}(U_{i})=U_{i}\times \prod_{j\neq i}X_{j}$. By Lemma \ref{12}, $\pi_{i}^{-1}(U_{i})\in ^{*}O(\prod_{i}X_{i})$.
\endproof

\begin{definition}
A map $f:X\to Y$ is called almost continuous if the inverse image of every regular open subset of $Y$ is open subset of $X$.
\end{definition}

\begin{remark}
It is clear that every $*-$continuous map is almost continuous. But, the converse need not be true. See the map $f$ of Example \ref{28}.
\end{remark}

\begin{lemma}
An open almost continuous map is $*-$continuous.
\end{lemma}

\proof
Let $f:X\to Y$ be an open almost continuous map and $y\in W\in ^{*}O(Y)$. Then, there exists $V\in R.O(Y)$ containing $y$ such that $V\subseteq W$. By \cite[Lemma 3.17]{M}, $f^{-1}(V)\in R.O(X)$ and $f^{-1}(V)\subseteq f^{-1}(W)$. So, $f^{-1}(W)\in ^{*}O(X)$.
\endproof










\begin{lemma}\label{29}
Let $f:X\to Y$ be a onto $*-$continuous map. If $X$ is connected, then $Y$ is connected.
\end{lemma}

\proof
Let $A$ be a clopen subset of $Y$. Then, $A$ is a clopen set in $(Y,^{*}O(Y))$. Hence, $f^{-1}(A)$ is a clopen set in $(X,^{*}O(X))$. By Lemma \ref{18}, $(X,^{*}O(X))$ is connected. So, $A=\emptyset$ or $A=Y$ and proof is complete.
\endproof

\begin{remark}
We know that compactness property preserved by continuous maps. But, this property do not preserve by $*-$continuous maps even though it is open. See Example \ref{30}.
\end{remark}

\begin{example}\label{30}
Let $(\mathbb{R},\tau)$ be as in Example \ref{2}. Consider the identity map $1_{\mathbb{R}}:(\mathbb{R},\tau_{1})\to (\mathbb{R},\tau)$. Then, $[0,1]$ is compact in $(\mathbb{R},\tau_{1})$, but it is not compact in $(\mathbb{R},\tau)$ ( Recall that the only compact sets in $(\mathbb{R},\tau)$ are finite).
\end{example}

\begin{remark}
Path connected and locally connected properties need not to be preserved by $*-$continuous maps. See Examples \ref{30},\ref{27} and \ref{24}.
\end{remark}

\begin{remark}
Recall that if $f:X\to Y$ is continuous, then $f(cl(A))\subseteq cl(f(A))$ for every subset $A$ of $X$. But, it is not true that $f(^{*}cl(A))\subseteq ^{*}cl(f(A))$. See Example \ref{31}.
\end{remark}

\begin{example}\label{31}
Consider the continuous map $f$ as in the Example \ref{28}. Then, $f(^{*}cl(\mathbb{Q}))=f(\mathbb{R})=(\mathbb{R-Q})\bigcup \{0\}$. But, $^{*}cl(f(\mathbb{Q}))=^{*}cl(\{0\})=\{0\}$.
\end{example}

\begin{lemma}
Let $f:X\to Y$ be a map. Then, $f$ is $*-$continuous if and only if for every subset $A$ of $X$, $f(^{*}cl(A))\subseteq ^{*}cl(f(A))$.
\end{lemma}

\proof
Let $x\in ^{*}cl(A)$ and $W$ be a regular open set in $Y$ contains $f(x)$. Then, $f^{-1}(W)$ is a $*-$open subset of $X$ contains $x$. Hence, there exists a regular open set $V$ in $X$ such that $f(x)\in f(V)\subseteq W$. Since $x\in ^{*}cl(A)$ , $V\bigcap A  \neq \emptyset $. So, $\emptyset \neq f(V)\bigcap f(A)\subseteq W\bigcap f(A)$. Hence, $x\in ^{*}cl(f(A))$. Conversely, let $W\in ^{*}O(Y)$. Then, $Y-W$ is $*-$closed. By assumption , $$f(^{*}cl(f^{-1}(Y-W)))=f(^{*}cl(X-f^{-1}(W))\subseteq ^{*}cl(Y-W)=Y-W$$
Hence, $^{*}cl(X-f^{-1}(W))\subseteq X-f^{-1}(W)$. So, $f^{-1}(W)\in ^{*}O(X)$ and the proof is complete.
\endproof

\end{document}